\documentclass[11pt]{amsart}
\usepackage{amsmath,amsfonts,latexsym,amssymb,amsthm, verbatim}
\usepackage{url}
 \usepackage[OT2, T1]{fontenc}

\usepackage{resizegather}
 \usepackage[flushleft]{threeparttable}
\usepackage{lscape}
\usepackage{amscd}

\usepackage{caption,booktabs}
\usepackage[]{amsrefs}
\usepackage{setspace}

\usepackage[flushleft]{threeparttable}

\usepackage{cleveref}

\newtheorem{thm}{Theorem}
\newtheorem{prop}[thm]{Proposition}
\newtheorem{lem}[thm]{Lemma}

\theoremstyle{definition}

\newcommand{\F}{\mathbb{F}}

\newcommand{\BQ}{\mathbb{Q}}

\newcommand{\BZ}{\mathbb{Z}}

\begin{document}

\title[]{Points of Order 13 on Elliptic Curves}

\author{S. Kamienny and B. Newman}
\address{\hspace{-5.8mm} Department of Mathematics\\ University of Southern California, \newline Los Angeles, CA 90089\\ USA}
\email{kamienny@usc.edu}
\email{bnewman235@gmail.com}



\maketitle

\section{Introduction}

	Our modern understanding of torsion points on elliptic curves over number fields began with a 1973 paper \cite{MT} of Mazur and Tate, from which we have borrowed the title of this note.  Spurred on by work of Ogg \cite{Ogg}, they carried out a descent in flat cohomology that proved that no elliptic curve over $\mathbb{Q}$ could possess a $\mathbb{Q}$-rational point of order 13.  A few years later Mazur gave his first proof of Ogg's Conjecture describing the possible rational torsion subgroups of elliptic curves over $\mathbb{Q}$.  Mazur's second proof \cite{Mazur} of Ogg's Conjecture lay the groundwork for an attack on the Strong Uniform Boundedness Conjecture that was eventually proved by Merel \cite{Merel}.  Merel's work provided an explicit bound on the size of the torsion subgroup over a degree $d$ number field, but it left open the problem of explaining the source of degree $d$ torsion when it does occur.  
    
	Here we return to the prime 13, the original subject of \cite{MT}.  The modular curve $X_1(13)$, whose non-cuspidal points classify isomorphism classes of elliptic curves with a rational torsion point of order 13, has genus 2 so there are infinitely many elliptic curves with a rational point of order 13 defined over quadratic number fields.  Bosman, Bruin, Dujella, and Najman  \cite{BBDN} have shown that any such curve must be defined over a real quadratic field.  Derickx, Kamienny and Mazur \cite{DKM} prove that any elliptic curve with a point of order 13 defined over a quadratic number field is part of a rationally parameterized family of such curves.  Unlike the original work of Mazur and Tate, where the use of an equation for $X_1(13)$ was expressly forbidden, the more recent work cited above has shown the value of using an equation.  We continue this approach here to study elliptic curves with points of order 13 defined over $\mathbb{Q}(\zeta_{13})^+$, and its quadratic extensions.
    

We show that no elliptic curve defined over $K=\mathbb{Q}(\zeta_{13})^+$ can possess a rational point of order 13.  We also find a finite number of elliptic curves with rational points of order 13 defined over quadratic extensions of $K$ that do not arise as part of a rationally parameterized family.  These curves owe their existence to the fact that $X_1(13)$ becomes bi-elliptic over $K$, and each of the elliptic curve factors contains finitely many $K$-rational points.  In a similar fashion $X_0(37)$ is a bi-elliptic curve that covers an elliptic curve with rank one over $\mathbb{Q}$.  This curve gives birth to an infinite family of elliptically parameterized elliptic curves with 37-isogenies defined over quadratic number fields, and this family is distinct from the infinite family that arises from the hyperellipticity of $X_0(37)$.

We are indebted to Filip Najman for suggesting to us that an earlier version of this
paper could be strengthened by eliminating the dependency on the weak Birch and
Swinnerton-Dyer conjecture. We would also like to thank Ken Ribet for clarifying
our understanding of modularity over totally real fields.

\section{The Modular Curve $X_1(13)$}

	The curve $X_1(13)$ is a cyclic cover of $X_0(13)$ with covering group $\Gamma$  isomorphic to $(\mathbb{Z}/13\mathbb{Z})^*/(\pm 1)$.   The automorphism group  of the curve $X_1(13)$ is a twisted dihedral group that is an extension of $\mathbb{Z}/2\mathbb{Z}$ by the group $\Gamma$.  Each of the involutions in the group  is a lift of the Atkin-Lehner involution of $X_0(13)$, and is defined over $K$.  Following Mazur and Tate \cite{MT} we denote these involutions by  $\tau_{\zeta}$, where  $\zeta$ is a primitive $13$th root of unity, and we identify the involutions associated with $\zeta$   and $\zeta^{-1}$.  The quotient of $X_1(13)$ by the action of  any  $\tau_{\zeta}$ is an elliptic curve defined over $K$.
    
	The curve $X_1(13)$ has twelve cusps,  six of them are $\mathbb{Q}$-rational, and the remaining 6 are rational over $K$.  The involutions  $\tau_{\zeta}$ interchange the two sets of cusps.  As is usual, we write $J_1(13)$ for the jacobian of  $X_1(13)$.  When we embed $X_1(13)$ in $J_1(13)$ the divisor classes supported at the rational cusps generate a $\mathbb{Q}$-rational subgroup $C$ of order 19.  The divisor classes supported at the six $K$-rational cusps generate a $K$-rational subgroup $D$, also of order 19.

The curve $X_1(13)$ has a model of the form $y^2=f(x)$ where 
$$f(x) = x^6 + 4x^5 + 6x^4 + 2x^3 + x^2 + 2x + 1$$
Magma tells us $ K= \mathbb{Q}(a)$  where $a$ is a root of 
$$x^6-x^5-5x^4+4x^3 + 6x^2-3x-1$$ 
Because $X_1(13)$ is bielliptic, it has a model of the form $$ y^2 = c_6x^6 + c_4x^4 + c_2x^2 + c_0$$
We will now describe how to obtain a model of the above form.
Let \\F=$\mathbb{\overline{Q}}[x,y]/(y^2-f)$ be the function field of $X_1(13)$  and let $e_0$ denote the hyperelliptic involution of $F$.  The fixed field of $e_0$ is generated by the image of $x$ in F.
An \textit{elliptic} involution of $G = $ Aut$(F)$ is an involution different
from $e_0$ \cite{Shaska}.  Using Magma to compute Aut($X_1(13)$) we found an involution $e$ which induces an automorphism of $F$ where
$$ x \mapsto \frac{x+b}{cx-1}  $$

$$b = -a^5 + 5a^3 - 6a$$
$$c =  -a^5 + 5a^3 - 6a - 1 $$
According to \cite{Shaska}, if one can find a a generator $X$ of the fixed field of $e_0$ such that $e(X) = -X$ then there is a relation of the form $$Y^2 = c_6X^6 + c_4X^4 + c_2X^2 + c_0$$ for some $Y \in F$.  We found that if $X$ is of the form:
$$ X = \frac{x+d_1}{x+d_2}   $$
with $d_1,d_2 \in \mathbb{\overline{Q}}$ then 
a simple calculation shows we will have $e(X) = -X$ if $$d_1 + d_2 = -2/c$$ $$d_1d_2  =  -b/c$$  Solving we obtain:
$$  d_1 =  -a^5 + 2a^4 + 3a^3 - 6a^2 + 1$$
$$ d_2 =  a^5 - 2a^4 - 5a^3 + 8a^2 + 6a - 5 $$
Using Magma to perform Gaussian elimination to find an appropriate linear combination of $X^6,X^4,X^2$ and $1$ we obtain the relation:
$$ Y^2 = c_6X^6 + c_4X^4 + c_2X^2 + c_0$$
where $$Y = \frac{y}{(x+d_2)^3} $$


$$ c_0 =  1/208(-52a^5 + 48a^4 + 240a^3 - 193a^2 - 218a + 175)  $$
$$ c_2 = 1/208(-18a^5 + 99a^3 + 38a^2 - 77a - 25)   $$
$$ c_4 = 
1/208(40a^5 - 238a^3 - 9a^2 + 296a + 69) $$
$$ c_6 = 
1/208(30a^5 - 48a^4 - 101a^3 + 164a^2 - a - 11)  $$

\noindent This model  has the  genus 1 quotient 
$$E: y^2 = c_6x^3 + c_4x^2 + c_2x + c_0$$
and twisting by a square in $K$ we obtain the model:
$$ E': y^2 = x^3 + bx^2 + cx + d$$

$$ b =  1/208(40a^5 - 238a^3 - 9*a^2 + 296a +69)  $$
$$ c  = 1/3328(-a^5 + 11a^4 + 38a^3 - 19a^2 - 73a - 16)  $$
$$ d = 1/692224(180a^5 - 4a^4 - 939a^3 - 30a^2 + 1092a + 252)   $$
Factoring the 19-division polynomial in Magma yields a point $P=(x,y)$ on $E'$ of order 19 with $$x = 1/208(134a^5 - 36a^4 - 681a^3 + 36a^2 + 799a + 181)$$
Magma tells us the (analytic) rank of both $E'(K)$ and  the other genus 1 quotient of $X_1(13)$ is 0.






\section{The Rank of $J_1(13)(K)$}
   
 	The jacobian $J_1(13)$ is irreducible over $\mathbb{Q}$. However, the bi-ellipticity of $X_1(13)$ over $K$ induces a splitting (up to isogeny) of $J_1(13)$ into the product of two elliptic curves over $K$.  These two curves are (up to isogeny) the two elliptic curve quotients of $X_1(13)$ over $K$.  In the language of \cite{KN} these elliptic curves are bi-elliptically linked, one arises from the other by the simple process described in \cite{KN}.
The two elliptic curve quotients of $K$ are modular (by \cite{Anni}), or as Ribet
has pointed out to us, simply because they are quotients of $J_1(13)$ over $K$. They each have rank 0 over $K$ by \cite{Yuan} (see also \cite{Gross}). It follows immediately that $J_1(13)$ also has rank zero over $K$.


\section{K-rational points on $X_1(13)$}

We now wish to determine $X_1(13)(K)$.  A search via Magma yields 12 points (which are in fact, the complete list of cusps on our model of $X_1(13)$).
As $X_1(13)$ is genus 2, Magma embeds it in (1,3,1)-weighted projective space, yielding two points at infinity $\infty_1 = (1,-1,0)$ and $\infty_2 =  (1,1,0)$.  
The cusps of $X_1(13)$ are $\infty_1$,$\infty_2$, $(0,\pm 1)$,$(-1,\pm 1)$ and\\

\noindent $(-a^5+4a^3+a^2-3a-1, \pm( -6a^5 - 6a^4 + 31a^3 + 19a^2 - 21a - 5))$      \\ 
        $(a^3 - a^2 - 3a + 2, \pm( -11a^5 + 18a^4 + 43a^3 - 66a^2 - 26a + 33))$      \\ 
    $(a^5 - 5a^3 + 6a , \pm(5a^5 - 4a^4 - 32a^3 + 5a^2 + 45a + 12))$  
\bigskip

    






\noindent Now $X_1(13)$ embeds into $ J := J_1(13)$ via the Abel-Jacobi map:
$$ X_1(13) \rightarrow J_1(13)$$
$$ P \mapsto P - \infty_1$$
\noindent Using the image of the cusps under the Abel-Jacobi map, we can generate $19^2$ torsion points in $J[19]$.
On the other hand, we will argue that $J(K)$ has size at most $19^2$.  
The discriminant of our model of $X_1(13)$ is $2^{12}\cdot 13^2$ and hence the curve remains nonsingular modulo primes of $\mathcal{O}_K$  above 3 and 5, so $J$ has good reduction at those primes.
Magma tells us 3 splits as the product of 2 distinct primes in $\mathcal{O}_K$ so 6=efg where g=2, e=1 and hence f = 3. Similarly Magma tells us 5 splits as the product of 3 distinct primes in $\mathcal{O}_K$ so  $g=3$, $e=1$ and hence $f = 2$. Magma tells us $|J(F_{27})| = 4\cdot 19^2 $ and $|J(F_{25})| = 19^2 $.

Since the reduction map  $J(K)_{tor} \rightarrow J(F_{25})$ (modulo a prime above 5) 
is injective on the prime-to-5 part of $J(K)_{tor}$, we see the prime-to-5 part of $J(K)_{tor}$ has size at most $19^2$.  On the other hand, under reduction mod 3 the 5-part of $J(K)_{tor}$ injects into $J(F_{27})$ and so must be trivial.



Now, as $X_1(13)$ is genus 2 (hence hyperelliptic), each element of $J_1(13)$ has a unique representative of the form $P + Q - (\infty_1 + \infty_2)$ where $P$ and $Q$ are points on $X_1(13)$ and if $P$ and $Q$ are affine then they don't lie on the same vertical line.  Each $K$-rational point $P$ on $X_1(13)$ gives rise to the following two points on $J_1(13)$:

$$P + \infty_1 - (\infty_1 + \infty_2)$$
$$P + \infty_2 - (\infty_1 + \infty_2)$$

As there are 12 points on $X_1(13)(K)$ this yields 23 = 24-1 points on $J_1(13)(K)$ as $\infty_1 + \infty_2 - (\infty_1 + \infty_2)=0_J$ is counted twice.  On the other hand, among the 361 elements of $J_1(13)(K)$, Magma tells us the Mumford representation $(a(x), b(x),d)$ satisfies deg(a(x))$<2$ for 23 elements (these are the elements in which $P$ or $Q$ is a point at infinity ). Hence the 12 points above are all the $K$-rational points on $X_1(13)$.

Each elliptic curve factor of $X_1(13)_{/K}$ has Mordell-Weil group (over $K$) isomorphic to $\mathbb{Z}/19\mathbb{Z}$. The inverse image of the 19 $K$-rational points on $E'$ gives us $38$ points on $X_1(13)_{/K}$ defined over quadratic extensions of $K$. Since $X_1(13)$ has only $12$ cusps we have found a collection of $26=38-12$ points of $X_1(13)_{/K}$ defined over quadratic extensions of $K$, but not defined over $K$. Each of these points corresponds to an elliptic curve with a point of order $13$ defined over a quadratic extension of $K$. Moreover, at least $25$ of these points must correspond to elliptic curves that do not arise as part of a rationally parameterized family (since, as Derickx has pointed
out to us, the image of $\mathbb{P}^1$ must map to a point in the jacobian $J_1(13)$).


\section{Quadratic Points on $X_0(37)$}

The curve $X_0(37)$ has a model of the form $y^2=f(x)$ where 
$$f(x) = (1/4)x^6+2x^5-5x^4 + 7x^3 - 6x^2 +3x - 1$$
Because $X_0(37)$ is bielliptic, it has a model of the form $$ y^2 = c_6x^6 + c_4x^4 + c_2x^2 + c_0$$

We will now describe how to obtain a model of the above form.
Let F=$\mathbb{\overline{Q}}[x,y]/(y^2-f)$ be the function field of $X_0(37)$  and let $e_0$ denote the hyperelliptic involution of $F$.  
Using Magma we found an involution $e$ which induces the automorphism of $F$ where
$$ x \mapsto \frac{x}{x-1}  $$
As mentioned before, if one can find a generator $X$ of the fixed field of $e_0$ such that $e(X) = -X$ then there is a relation of the form $$Y^2 = c_6X^6 + c_4X^4 + c_2X^2 + c_0$$ for some $Y \in F$.  We found 
$$ X = \frac{x-2}{x}   $$

\noindent satisfies $e(X) = -X$.Using Magma 
we obtain the relation:
$$ Y^2 = (-1/64)(X^6 + 9X^4 + 11X^2 -37)$$
where $$Y = \frac{y}{x^3} $$
This model  has the obvious  genus 1 quotient 
$$ Y^2 = (-1/64)(X^3 + 9X^2 + 11X -37)$$
This curve is isomorphic over $\mathbb{Q}$ to the elliptic curve:
$$ E: y^2 + y = x^3 - x   $$
and Magma tells us $E(\mathbb{Q})$
is rank 1 with trivial torsion subgroup over $\mathbb{Q}$.  Furthermore Magma tells us $E(\mathbb{Q}) = <(0,0)>$.  %
Using Magma we obtain the quotient map from Magma's model of $X_0(37)$ to the curve $E$:
    $$ (x,y) \mapsto (\frac{-x^2+x-1}{x^2},\frac{-x^3+y}{x^3}   )  $$
The inverse image of the Mordell-Weil group $E(\mathbb{Q})$ under this map gives us an
infinite number of quadratic points on $X_0(37)$, and as before, at most one of these points can come from a rationally parameterized family.  Using Magma, we wrote a program which generates these quadratic points and for each such point
finds an elliptic curve with a 37-isogeny in the corresponding isomorphism class\footnote{https://github.com/bdnewman/phd-projects} (except for isomorphism classes corresponding to elliptic curves with $j = 0$ or $1728$).  \Cref{37} below lists some examples of quadratic points on $X_0(37)$.
Some of the curves in \Cref{37} already appear in the tables of Cremona \cite{LMFDB}.

\begin{table}[ht]
\caption{ Elliptic curves with 37-isogenies defined over\\ quadratic fields}


\label{37}

\renewcommand{\arraystretch}{1.7}

\hspace*{-4cm}
 \begin{tabular}{| l | l | l | }
    \hline
    $D$ &  $P$ & $E$ \\ \hline
    -3 & $[\frac{(-a + 1)}{2},0]$ & $[\frac{315a - 3285}{2},
    3630a - 24948]$   \\ \hline
    
     -7 & $[\frac{(-a + 1)}{4},0]$ &  $[\frac{-765a + 6345}{8},
    \frac{-30753a - 40635}{8}]$
       \\ \hline
 
  -11 & $[\frac{-a + 1}{6},\frac{a - 4}{9}]$ & $[\frac{-640a + 2848}{3},
    \frac{75040a + 356048}{27}]$ 
       \\ \hline
       
       -1 & $[\frac{-4a + 2}{5},\frac{2a - 11}{25}]$ & $j = 1728$  
       \\ \hline
       
       -3 & $[\frac{-3a + 1}{14}, \frac{-18a + 20}{49}]$ & $j = 0$ 
       \\ \hline
       
       -7 &  $[\frac{-3a + 9}{8},\frac{5a + 9}{16}]$ & $[ \frac{34425a - 207315}{32},
   \frac{3224205a - 11925711}{64}]$  
       \\ \hline
       
 		-159 & $[\frac{-5a + 25}{92},\frac{63a - 1695}{8464}]$&  $\frac{1992776643585a + 394997768625969}{274877906944}$ 
       \\ \hline
       
       -67 & $[\frac{-7a + 49}{58},\frac{-150a - 1995}{841}]$ & $[\frac{126720a - 16964640}{841},
    \frac{-301386960a + 26856906048}{24389}]$
       \\ \hline
       
       -173 & $[\frac{-4a + 8}{177},\frac{5635a - 36050}{93987}]$
&

$[\frac{-2366000*a + 214423015}{93987},
    \frac{473705169712*a + 1346165714530}{149721291}]$  
       \\ \hline
 
 	     -2051 &  $[\frac{-23a + 529}{1290},\frac{1624a - 452732}{2080125}]$ & 
         
         $[\frac{-269113a - 6547514791}{6933750},
    \frac{-3221900558162a - 383038176258584}{33542015625}]$
       \\ \hline

  -7951 & $[\frac{-29a + 841}{4396}, \frac{-31211a + 3837251}{16909214}]$ &

$[ \frac{2989314135a + 825198488937}{473457992},
    \frac{-1023161515376435a + 46768183795198699}{3642312332456}   ] $
     \\ \hline

      4521 & $[\frac{-59a - 3481}{520},\frac{-2495247a - 167683133}{540800}]$&
     
$[\frac{-12322582501a - 1272066914239}{10816000},
    \frac{20998975870933249a + 1250157035949620211}{50618880000}]$     
      
      \\ \hline 
      
       -124027& $[\frac{-129a + 16641}{70334},
    \frac{15848993a - 9806545170}{28444511447}]$ &
    
    $[\frac{-3683948697936a + 792633701552976}{2081621064985},
    \frac{339400774163819886912a + 625337457592756853499120}{92603525510294281175}]$
    
    \\ \hline


    \end{tabular}
\hspace*{-3cm}




$\dagger$ For each quadratic point $P$ on $X_0(37)$ (defined over $\mathbb{Q}(a)$ with $a^2 = D$) listed above, we provide 
the coefficients 
$[A,B]$ of an elliptic curve with model $ y^2=x^3 + Ax+B$ in the isomorphism class corresponding to $P$.
\end{table}

\newpage

\section{Appendix:  Points of Order $13$ on Elliptic curves over $\mathbb{Q}(\zeta_{13})^+$ --  A second approach}

 Let $K=\mathbb{Q}(\zeta_{13})^{+}$, $\mathcal{O}$ its ring of integers, and suppose that $x=(E,P)$ is a $K$-rational point on $X_1(13)$. We write $J$ for the Neron model of $J_1(13)$ over $S$ = Spec $\mathcal{O}$. We write $\mathfrak{p}$ for a prime of $K$, $k(\mathfrak{p})$ for its residue field, and $p$ for the characteristic of $k(\mathfrak{p})$. In the following we work over the base $S$.  
 
 If $E$ has additive reduction at $\mathfrak{p}$ then the point $P$ (of order $13$) must reduce to $(E/k(\mathfrak{p}))^o$, since $[E:E^o]$ is bounded by $4$.  However, $(E/k(\mathfrak{p}))^o$  is an additive group, and an additive group in characteristic $p$ cannot have a point of order $13$ unless $p = 13$. 
 
 If $E$ has (potentially) multiplicative reduction at $\mathfrak{p}$ then $x$ must reduce (mod $\mathfrak{p}$) to a cusp $Q$   of $X_1(13)$.  The class of $(x -Q)$ is a $K$-rational point $T$ on $J_1(13)$, and hence is torsion. The point $T$ generates a finite flat subgroup scheme $C$ of $J$.  Since $J_1(13)(K)$ has no 2-torsion this group scheme must be {\'e}tale.  The point  $x$ reduces (mod $\mathfrak{p}$) to the cusp $Q$, so the group scheme $C$ reduces (mod $\mathfrak{p}$)  to 0.  Because C is {\'e}tale it must already be 0 in characteristic 0, i.e. $(x-Q)$ is the divisor of a function on the genus 2 curve $X_1(13)$.  This is clearly impossible.  
 
 Thus, $E$ has good reduction at all primes of residue characteristic different from $13$, and at the prime above $13$ $E$ has good reduction or potentially good reduction.  If $E$ has good reduction at $13$ then it cannot exist by \cite{Yasuda}.  If $E$ has potentially good reduction then we use Cremona's algorithm (which we could also use in the case of good reduction) to find the finite set of possible $E$, and check that none of them possess a $K$-rational  point of order 13.  Of course, in this case we already know that none of the curves that we find will possess a $K$-rational point of order 13.  

	The advantage of this method is that it doesn't depend upon the 
hyperellipticity of $X_1(13)$, and it  will either prove that there are no elliptic curves over $K$ with a point of order 13, or it will find all such curves if they happened to exist.  The difficulty with this method is that it requires finding all integral points on certain associated elliptic curves, and this may be computationally (but not 
theoretically) prohibitive.

\end{document}